\documentclass{article} 
\usepackage{graphicx}
\usepackage{amsmath}
\usepackage{amsfonts}
\usepackage{amssymb}
\usepackage{latexsym}
\usepackage{tikz}
\usepackage{fancyvrb}
\usepackage{datetime}

\usepackage{graphicx} 
\usepackage{color}
\usepackage[colorlinks]{hyperref}

\newtheorem{lemma}{Lemma}[section]

\newtheorem{corollary}{Corollary}[section]
\newtheorem{remark}{Remark}
\newtheorem{theorem}{Theorem}[section]

\usepackage{graphicx}
\usepackage{amsmath}
\usepackage{amsfonts}
\usepackage{amssymb}
\usepackage{latexsym}
\usepackage{cite}
\begin{document}

	\begin{center}
{\Large{Graphs determined by signless Laplacian Spectra}}

\medskip	

{Ali Zeydi Abdian\footnote{Lorestan University, College of Science, Lorestan, Khoramabad, Iran; e-mail: abdian.al@fs.lu.ac.ir;  aabdian67@gmail.com; azeydiabdi@gmail.com},  Afshin Behmaram
\footnote{Faculty of Mathematical Sciences, University of Tabriz, Tabriz, Iran; Email: behmaram@tabrizu.ac.ir} and Gholam Hossein Fath-Tabar\footnote{Department of Pure Mathematics, Faculty of Mathematical Sciences, University of Kashan, Kashan 87317-53153, Iran; Email: fathtabar@kashanu.ac.ir.}
	
}	\end{center}

\medskip 
\begin{abstract} 
In the past decades, graphs that are determined by their spectrum have received more attention, since they have been applied to several fields, such as randomized algorithms, combinatorial optimization problems and machine learning. An important part of spectral graph theory is devoted to determining whether given graphs or classes of graphs are determined by their spectra or not. So, finding and introducing any class of graphs which are determined by their spectra can be an interesting and important problem. A graph is said to be $DQS$ if there is no other non-isomorphic graph with the same signless Laplacian spectrum. For a $DQS$ graph $G$, we show that $ G\cup rK_1\cup sK_2 $ is $DQS$ under certain conditions, where $ r $, $ s $ are natural numbers and $ K_1 $ and $ K_2 $ denote the complete graphs on one vertex and two vertices, respectively. Applying
these results, some $DQS$ graphs with independent edges and isolated vertices are obtained.\\
\textbf{Keywords:} Spectral characterization; Signless Laplacian spectrum; Cospectral graph.\\ 
\textbf{MSC(2010):} 05C50.
\end{abstract}

\section{Introduction}
Let $G = (V,E)$ be a simple graph with vertex set $ V = V (G) = \left\{ {v_1, . . . , v_n} \right\}$ and edge set $E = E(G) = \left\{ {e_1, . . . , e_m} \right\}$. Denote by $d(v)$ the degree of vertex $ v $. All graphs considered here are simple and undirected. All notions on graphs that are not defined here can be found in \cite{LP, Ba, B, CRS, W}. The join of two graphs $ G $ and $ H $ is a graph formed from disjoint copies of $ G $ and $ H $ by connecting each vertex of $ G $ to each vertex of $ H $. We denote the join of two graphs $ G $ and $ H $ by $ G\bigtriangledown H $. The complement of a graph $ G $ is denoted by $ \overline{G} $.\\

Let $A(G)$ be the $(0,1)$-adjacency matrix of graph $ G $. The characteristic polynomial of $G$ is $\det (\lambda I - A(G))$, and it is denoted by $ P_G(\lambda) $. Let $ \lambda_1 $, $ \lambda_2 $, ..., $\lambda_n$ be the distinct eigenvalues of $G$ with multiplicities $m_1$, $ m_2$, ... , $m_n$, respectively. The multi-set of eigenvalues of $ Q(G) $ is called the signless Laplacian spectrum of $G$. The matrices $ L(G) =D(G)-A(G) $ and $ Q(G)=SL(G) = D(G)+A(G) $ are called the Laplacian matrix and the signless Laplacian matrix of $ G $, respectively, where $D(G)$ denotes the degree matrix. Note that $D(G)$ is diagonal. The multi-set $ {\rm{Spec}}_{Q}(G) = \left\{ { [\lambda_1]^{ m_1 } , [\lambda_2]^{ m_2 } , ..., [\lambda _n]^{ m_n } } \right\}$ of eigenvalues of $ Q(G) $ is called the signless Laplacian spectrum of $G$, where $ m_i $ denote the multiplicities of $ \lambda_i $. The Laplacian spectrum is defined analogously.\\

 For any bipartite graph, its $Q$-spectrum
coincides with its $L$-spectrum. Two graphs are $Q$-cospectral (resp. $L$-cospectral, $A$-cospectral) if they have the
same $Q$-spectrum (resp. $L$-spectrum, $A$-spectrum). A graph $G$ is said to be $DQS$ (resp. $DLS$, $DAS$) if there is
no other non-isomorphic graph $Q$-cospectral (resp. $L$-cospectral, $A$-cospectral) with $G$. Van Dam and Haemers \cite{VH} conjectured that almost all graphs are determined by their spectra. Nevertheless, the set of graphs that are known to be detrmined by their spectra is too small. So, discovering infinite classes of graphs that are determined by their spectra can be an interesting problem. About the background of the question "Which graphs are determined by their spectrum?", we refer to \cite{VH}. It is interesting to construct new $DQS$ ($DLS$) graphs from known $DQS$ ($DLS$) graphs. For a $DLS$ graph $G$, the join $G\cup rK_1$ is also $DLS$ under some conditions \cite{LC}. Actually, a graph is $DLS$ if
and only if its complement is $DLS$. Hence we can obtain $DLS$ graphs from known $DLS$ graphs by adding  independent edges.

 \vspace{2mm} 
 
 Up to now, only some graphs with special structures are shown to be  {\it determined by their spectra} (DS, for short) (see \cite{A, AA, AAA, AAAA, AAAAA, AAAAAA, AAAAAA1, AAAAAA2, AAAAAA3,CHA, CRS1, CS, Das, HLZ, LP,  Mer, M, M1,SZ, WYY,  W1} and the references cited in them). \vspace{1mm}
\\

  In this paper, we investigate signless Laplacian spectral characterization of graphs with independent edges and isolated vertices. For a $DQS$ graph $G$, we show that $G \cup rK_1\cup sK_2$ is $DQS$ under certain conditions. Applying these results, some $DQS$ graphs with independent edges and isolated vertices are obtained.\\

\section{Some definitions and preliminaries}

Some useful established results about the spectrum are presented in this section,
 will play an important role throughout this paper.

\begin{lemma}\cite{A,  CHA, CRS}\label{lem 2-1} For the adjacency matrix, the Laplacian matrix and the signless Laplacian matrix of a graph $G$, the
following can be deduced from the spectrum:\\

(1) The number of vertices.

(2) The number of edges.

(3) Whether G is regular.

For the Laplacian matrix, the following follows from the spectrum:

(4) The number of components.

For the signless Laplacian matrix, the following follow from the spectrum:

(5) The number of bipartite components.

(6) The sum of the squares of degrees of vertices.
\end{lemma}

\begin{lemma}\cite{CHA}\label{lem 2-20} Let $ G $ be a graph with $ n $ vertices, $ m $ edges and $ t $ triangles and vertex degrees $ d_1, d_2, . . . , d_n $. Let $ {T_k} = \sum\limits_{i = 1}^n {{{({q_i}(G))}^k}}$, then

$ T_0=n $, $ T_1=\sum\limits_{i = 1}^n d_i=2m $, $ T_2=2m+\sum\limits_{i = 1}^n d^2_i $ and $ T_3=6t+3\sum\limits_{i = 1}^n d^2_i+\sum\limits_{i = 1}^n d^3_i $.
\end{lemma}

For a graph $G$, let $P_L(G)$ and $P_Q(G)$ denote the product of all nonzero eigenvalues of $L_G$ and $Q_G$, respectively. We assume that $P_L(G) = P_Q(G) = 1$ if $G$ has no edges. 

\begin{lemma}\cite{CRS}\label{lem 2-3}  For any connected bipartite graph $G$ of order $n$, we have $P_Q(G) = P_L(G) = n\tau(G)$, where $\tau(G)$ is the
number of spanning trees of $G$. 
\end{lemma}
For a connected graph $G$ with $n$ vertices and $m$ edges, $G$ is called unicyclic (resp. bicyclic) if $m = n$ (resp.
$m = n + 1$). If $G$ is a unicyclic graph that contains  an odd (resp. even) cycle, then $G$ is called odd unicyclic (resp. even unicyclic).

\begin{lemma}\cite{MK}\label{lem 2-4} For any graph $G$, $det(Q_G) = 4$ if and only if $G$ is an odd unicyclic graph. If $G$ is a non-bipartite connected graph and $|E(G)|>|V(G)|$, then $det(Q_G) > 16$, with equality if and only if $ G $ is a non-bipartite bicyclic graph with $ C_4 $ as its induced subgraph.
\end{lemma}

\begin{lemma}\cite{CZ}\label{lem 2-5} Let $H$ be a proper subgraph of a
connected graph $G$. Then, $q_1(G) > q_1(H)$.
\end{lemma}
\section{Main Results} 

We first investigate spectral characterizations of the union of a tree and several complete graphs $ K_1 $ and $ K_2 $.
\begin{theorem}\label{the 3-1}
 Let $T$ be a $DLS$($DQS$) tree of order $n$. Then $T\cup rK_1\cup sK_2$ is $DLS$.\\
 
 $T\cup rK_1\cup sK_2$ is $DQS$  if $n$ is not divisible by 2 and $ s=1 $.
 \end{theorem}
\noindent \textbf{Proof} Let $G$ be any graph $L$-cospectral with $T\cup rK_1\cup sK_2
$. By Lemma \ref {lem 2-1}, $G$ has $n +r+2s$ vertices, $n-1+s$ edges and
$r +s+ 1$ components. So each component of $G$ is a tree. Suppose that $G = G_0 \cup G_1 \cup . . . \cup  G_{r+s}$, where $G_i$ is a
tree with $n_i$ vertices and $n_0 \geq n_1 \geq . . . \geq n_s\geq ... \geq n_{r+s} \geq 1$. Since $G$ is $L$-cospectral with $T\cup rK_1\cup sK_2$, by Lemma \ref {lem 2-3}, we
get $n_0n_1... n_{r+s} = P_L(G) = n2^s$. We claim that $n_s=2$. Suppose not and so $n_s\geq 3$. Therefore, $n_0 \geq n_1 \geq . . . \geq n_s\geq 3$ and since $n_{s+1} \geq ... \geq n_{r+s}\geq 1$, one may deduce that $n2^s=n_0n_1... n_{r+s}\geq 3^{s+1}$ or $n(\dfrac{2}{3})^s\geq 3$. Now if $s\longrightarrow \infty$, then $0\geq 3$, a contradiction. Hence $n_s=2$. By a similar argument one may show that $n_1=n_2=...=n_{s-1}=2$ and so $n_0=n$ and $n_{s+1}=n_{s+2}=...=n_{s+r}=1$. Hence $G = G_0 \cup rK_1\cup sK_2$. Since $G$ and $T\cup rK_1\cup sK_2$ are $L$-cospectral, $G_0$ and $T$ are
$L$-cospectral. Since $T$ is $DLS$, we have $G_0 = T$, $G = T\cup rK_1\cup sK_2$. Hence $T\cup rK_1\cup sK_2$ is $DLS$.
Let $H$ be any graph $Q$-cospectral with $T\cup rK_1\cup sK_2$. By Lemma \ref{lem 2-1}, $H$ has $n +r+2s$ vertices, $n -1+s$ edges and $r+s+ 1$
bipartite components. So one of the following holds:\\

(i) $H$ has exactly $r+s+ 1$ components, and each component of $H$ is a tree.\\

(ii) $H$ has $r+s+ 1$ components which are trees, the other components of $H$ are odd unicyclic.\\

If (i) holds, then $H$ and $T\cup rK_1\cup sK_2$ are both bipartite, so they are also $L$-cospectral. Since $T\cup rK_1\cup sK_2$ is $DLS$, we
have $H = T\cup rK_1\cup sK_2$.\\ 

If (ii) holds, then by Lemma \ref {lem 2-4}, $P_Q(H)$ is divisible by 4. Since $T$ is a tree of order $n$, by
Lemma \ref {lem 2-3}, $P_Q(H)= n2^s$ is divisible by 4. Hence $T\cup rK_1\cup sK_2$ is $DQS$ when $n$ is not divisible by 2 and $ s=1 $.$ \Box $\\

\begin{remark} Some $DLS$ trees are given in \cite{Aa, Bo, Bu, Lu, Shen, Stani}. We can obtain $DLS$ ($DQS$) graphs with independent
edges and isolated vertices from Theorem \ref{the 3-1}.
\end{remark}

\begin{theorem}\label{the 3-2}
 Let $G$ be a $DQS$ odd unicyclic graph of order $n\geq 7$. Then $G \cup rK_1\cup sK_2$ is $DQS$.
 \end{theorem}
\noindent \textbf{Proof}  Let $H$
be any graph $Q$-cospectral with $G \cup rK_1\cup sK_2$. By Lemma \ref {lem 2-4}, $P_Q(H) = 4(2^s)$. By Lemma \ref {lem 2-1}, $H$ has $n + r+2s$
vertices, $n+r$ edges and $r+s$ bipartite components. So one of the following holds:\\

(i) $H$ has exactly $r+s$ components, and each component of $H$ is a tree.\\

(ii) $H$ has $r+s$ components which are trees, the other components of $H$ are odd unicyclic.\\

If (i) holds, then we can let $H = H_1 \cup . . . \cup H_{r+s}$, where $H_i$ is a tree with $n_i$ vertices and $n_1 \geq . . . \geq n_{r+s} \geq 1$. Since $P_Q(H)= 4(2^s)$, by Lemma \ref {lem 2-3}, we have $n_1. . .n_{r+s} = 4(2^s)$, $n_1\leq 8 $. \\

Since $G$ contains a cycle, we have $q_1(H) = q_1(G) \geq 4$. Let $\Delta(H)$ be the maximum degree of $H$. If $\Delta(H)\leq 2$, then all components of $H$ are paths, i.e., $q_1(H) < 4$, a contradiction. So $\Delta(H)> 3$. From $n_1\leq 8$ and $n_1 . . . n_{r+s} = 4(2^s)=2^{(s+2)}$, we know that $H_1 = K_{1,7}$ (without loss of generality),
$H_2 = . . . = H_s = K_2$ and $H_{s+1} = . . . = H_{r+s} = K_1$. Since $H = K_{1,7} \cup (s - 1)K_2\cup rK_1$ has $n + r+2s$ vertices, we get $n = 6$, a contradiction to $n > 6$.\\

If (ii) holds, then we can let $H = U_1\cup . . . \cup Uc \cup H_1\cup . . . \cup H_r$, where $U_i$ is odd unicyclic, $H_i$ is a tree with
$n_i$ vertices. By Lemma \ref {lem 2-3} and \ref{lem 2-4}, $4(2^s)= P_Q(H) = 4^cn_1. . . n_r$. So $c = 1$, $H_1 = . . . = H_s= K_2$ and $H_{s+1} = . . . = H_{r+s}= K_1$.  Since $H = U_1\cup rK_1\cup sK_2$ and $G \cup rK_1\cup sK_2$ are $Q$-cospectral, $U_1$ and $G$ are $Q$-cospectral. Since $G$ is $DQS$, we have $U_1 = G$, $H = G \cup rK_1\cup sK_2$.$ \Box $

\begin{remark} Some $DQS$ unicyclic graphs are given in \cite{Bu1, Liu2, Liu3, Liu4, WAN, Zhang}. We can obtain $DQS$ graphs with independent edges and isolated vertices from Theorem \ref{the 3-2}.
\end{remark}

\begin{theorem}\label{the 3-3}
Let $G$ be a non-bipartite $DQS$ bicyclic graph with $C_4$ as its induced subgraph and $ n\geq 5 $. Then $ G\cup rK_1\cup sK_2 $ is $DQS$.
\end{theorem}
\noindent \textbf{Proof} Let $H$ be any graph $Q$-cospectral with $G\cup rK_1\cup sK_2$. By Lemma \ref {lem 2-4}, we have $P_Q(H) = 16(2^r)$. By
Lemma \ref {lem 2-1}, $H$ has $n +r+ 2s$ vertices, $n + 1+s$ edges and $r+s$ bipartite components, where $n = |V(G)|$. So $H$ has at least
$r+s-1$ components which are trees.
Suppose that $H_1,H_2, . . . , H_{r+s}$ are $r+s$ bipartite components of $H$, where $H_2, ... , H_r$ are trees. If $H_1$ contains an
even cycle, then by Lemma \ref {lem 2-3}, we have $P_Q(H) \geq P_Q(H_1) \geq 16$, and $P_Q(H) = 16(2^{s-1})$ if and only if $H = C_4\cup (s-1)K_2\cup rK_1$. Since $H$ has $n + r+2s$ vertices, we get $n = 2$, a contradiction ($G$
contains $C_4$). Hence $H_1,H_2, . . . , H_{r+s}$ are trees.
Since $H$ has $n +r+2s$ vertices, $n + 1+s$ edges and $r+s$ bipartite components, $H$ has a non-bipartite component
$H_0$ which is a bicyclic graph. Lemma \ref {lem 2-4} implies that $P_Q(H) > P_Q(H_0) > 16$, and $P_Q(H) = 16(2^s)$ if and only if
$H = H_0 \cup rK_1\cup sK_2$ and $H_0$ contains $C_4$ as its induced subgraph. By $PQ(H) = 16(2^s)$, we have $H = H_0 \cup rK_1\cup sK_2$. Since $H$
and $G \cup rK_1\cup sK_2$ are $Q$-cospectral, $H_0$ and $G$ are $Q$-cospectral. Since $G$ is $DQS$, we have $H_0 = G, H = G \cup rK_1\cup sK_2$.
Hence $G \cup rK_1\cup sK_2$ is $DQS$.$ \Box $

\begin{remark}Some $DQS$ bicyclic graphs are given in \cite{Guo, Liu5, WAN1, WAN2}. We can obtain $DQS$ graphs with independent
edges and isolated vertices from Theorem \ref{the 3-3}.
\end{remark}

\begin{theorem}\label{the 3-4}
Let $G$ be a $DQS$ {\bf connected} non-bipartite graph with $n\geq 3$ vertices.   If $H$ is $Q$-cospectral with $G \cup rK_1\cup sK_2$, then $ H $ is a $DQS$ graph.\\
\end{theorem}

\noindent \textbf{Proof} By Lemma \ref {lem 2-1}, $H$ has $n + r+2s$ vertices and at least $r+s$ bipartite components. We perform the mathematical induction on $ s $.\\

 $H$ has $r+s$ components. Since $H$ has at least $r+s$ bipartite components, each component of $H$ is bipartite. Suppose that $H = H_1\cup ... \cup H_{r+s}$, where $H_i$ is a connected bipartite graph with $n_i$ vertices, and
$ n_1 \geq . . . \geq n_s\geq ... \geq n_{r+s} \geq 1$. Since $H$ and $G\cup rK_1\cup sK_2$ are $Q$-cospectral, by Lemma \ref {lem 2-1}, $G$ is a connected non-bipartite graph.\\

Let $ s=1 $. For $ n\geq 3 $, $ q_1(G)\geq 3$, since $ G $ has $ K_{1,2} $ or $ K_3 $ as its subgraph.  Obviously $ {\rm{Spec}}_{Q}(H)$ has exactly $ r+s $ eigenvalues that are zero. We show that if $H$ is $Q$-cospectral with $G \cup rK_1\cup K_2$, then $ H $ is a $DQS$ graph. First we show that there is no connected graph $Q$-cospectral with ${\rm{Spec}}_{Q}(G^{'})={\rm{Spec}}_{Q}(G) \cup \left\{ {{{\left[ 2 \right]}^1}} \right\}$. In fact we prove that $ G^{'} $ cannot have 2 as its eigenvalue. Obviously, ${\rm{Spec}}_{Q}(H)={\rm{Spec}}_{Q}(G^{'}) \cup \left\{ {{{\left[ 0 \right]}^{r+1}}} \right\}$. But, in this case $ |E(G^{'})|=|E(G)|+1$ and $ |V(G^{'})|=|V(G)|+1$, which means that $G^{'}$ must be  connected. Otherwise, $G^{'}$ contains 0 as its signless eigenvalues, a contradiction. Therefore, $ G $ is a proper subgraph of $G^{'}$ and  so $ q_1(G^{'})\gneqq q_1(G)\geq 3$ (see Lemma \ref{lem 2-5}), a contradiction. Therefore, $ G^{'} $ cannot have 2 as its eigenvalue. 
By what was proved one can easily conclude that ${\rm{Spec}}_{Q}(H)={\rm{Spec}}_{Q}(G) \cup {\rm{Spec}}_{Q}(K_2)\cup {\rm{Spec}}_{Q}(rK_1)$, since $ G $ is not a bipartite graph and so has not $ 0 $ as an its signless Laplacian eigenvalue. Therefore, $ H=G\cup K_2\cup rK_1 $.\\ 

Now, let the theorem be true for $ s $; that is, if ${\rm{Spec}}_{Q}(G_1)={\rm{Spec}}_{Q}(G) \cup {\rm{Spec}}_{Q}(rK_1\cup sK_2)$, then $ G_1=G\cup rK_1 \cup sK_2$. We show that it follows from ${\rm{Spec}}_{Q}(K)={\rm{Spec}}_{Q}(G) \cup {\rm{Spec}}_{Q}(rK_1\cup (s+1)K_2)$ that $ K=G\cup rK_1\cup (s+1)K_2
$. Obviously, $ K $ has $ 2 $ vertices, one edge and one bipartite component more than $ G_1 $. So, we must have $ K=G_1\cup K_2 $. $ \Box $\\

 \begin{remark}
 In the following results graph $ G $ in $ G\cup rK_1\cup sK_2 $ is a connected non-bipartite.
 \end{remark}
 \begin{corollary}The graph $K_n\cup rK_1\cup sK_2$ is $DQS$.
 \end{corollary}
 \noindent \textbf{Proof} From \cite{VH} (Proposition 7), if $ n=1, 2 $, then $ K_n\cup rK_1\cup sK_2 $ is $ DQS $. For $ n\geq 3 $, by Theorem \ref{the 3-4} the result follows.
 $ \Box $
 \\
 
 In \cite{CH}, C\'amara and Haemers proved that a graph obtained from $ K_n $ by deleting a matching is $DAS$. In \cite{FIL}, it have been shown that this graph is also $DQS$.
 \\
 
 \begin{corollary}\label{cor 3-2}
 Let $ G $ be the graph obtained from $ K_n $ by deleting a matching. Then $ G\cup rK_1\cup sK_2 $ is $DQS$. 	
 \end{corollary}
  \noindent \textbf{Proof} From \cite{VH} (Proposition 7), if $ n=1, 2 $, then $ K_n\cup rK_1\cup sK_2 $ is $ DQS $. For $ n\geq 3 $, by Theorem \ref{the 3-4} the result follows.
 $ \Box $
 
 A regular graph is $DQS$ if and only if it is $DAS$ \cite{VH}. It is known that a $k$-regular graph of order $n$ is $DAS$ when $k = 0, 1, 2, n-1, n-2, n-3$ \cite{CHA}. Hence a $k$-regular graph of order $n$ is $DQS$ when $k = 0, 1, 2, n-1, n-2, n-3$.
 
 \begin{corollary}
 Let $ G$ be a connected $(n-2)$-regular graph of order $ n $. Then $ G\cup rK_1\cup sK_2 $ is $DQS$.
 \end{corollary}
 
 \begin{corollary}\label{cor 2}
 Let $ G$ be a connected $(n-3)$-regular graph of order $ n $. Then $ G\cup rK_1\cup sK_2 $ is $DQS$.
 \end{corollary}
 
\begin{corollary}
 Let $ G$ be a connected $(n-4)$-regular $DAS$ graph. Then $ G\cup rK_1\cup sK_2 $ is $DQS$.
 \end{corollary}

\begin{remark} Some $3$-regular $DAS$ graphs are given in \cite{VH, Liu6}. We can obtain $DQS$ graphs with independent edges and isolated vertices and isolated vertices from Corollary \ref{cor 2}.
\end{remark}

\begin{corollary} Let $ F_n $ denotes the friendship graph and $ G $ be $Q$-spectral with $ F_n $, then $ G\cup rK_1\cup sK_2$ is $ DQS $.
\end{corollary}
\noindent \textbf{Proof} It is well-known that $ F_n $ is $ DQS $. By Theorem \ref{the 3-4} the proof is completed.$ \Box $

\today
\end{document}